\documentclass[11pt, review]{elsarticle}

\usepackage{lineno,hyperref}
\modulolinenumbers[1]

\journal{XXXXXXXX}

\usepackage{geometry}
\usepackage{graphicx}
\usepackage{xcolor}

\usepackage{amsmath,amsthm}
\usepackage{mathrsfs}
\usepackage{amssymb,amsbsy} 
\usepackage{pst-all}
\psset{dotscale=1.2 1.2}

\makeatletter
\renewenvironment{proof}[1][\proofname]{\par
	\normalfont
	\topsep6\p@\@plus6\p@ \trivlist
	\item[\hskip4.35\labelsep\itshape
	#1\@addpunct{.}]\ignorespaces
}{%
	\qed\endtrivlist
}
\makeatother

\usepackage{multirow}
\usepackage{array}
\usepackage{booktabs}
\usepackage{dsfont}
\usepackage{enumerate}

\usepackage{blkarray}

\usepackage[mathscr]{euscript}

\def\bpsp{\begin{pspicture}}
\def\epsp{\end{pspicture}}

\renewcommand\c[1]{{\mathcal{#1}}}
\newcommand\matr[2]{{\left[\ba{#1}#2\ea\right]}}

\newcommand\orient[1]{\scalebox{1.8}{\rotatebox{#1}{\psline{->}(0,0)(.1,0)}}} 
\newcommand\mysubscript[1]{{\raisebox{-.1cm}{$#1$}}}
\newcommand\sdv[1]{\boldsymbol{\mu}_\mysubscript{#1}}
\newcommand\qsdv[1]{{}^{\tt q}\kern-2.5pt\boldsymbol{\mu}_\mysubscript{#1}}

\DeclareMathOperator*{\bdf}{\tt bd}
\DeclareMathOperator*{\bdq}{\tt bd_q}

\DeclareMathOperator*{\dist}{\tt dist}

\DeclareMathOperator*\diff{\ts \tt diff}

\def\bbm{\begin{bmatrix}}
	\def\ebm{\end{bmatrix}}

\newtheorem{thn}{Theorem}
\newcommand\bt{\begin{thn}}
	\def\et{\end{thn}}
\newtheorem{lemma}[thn]{Lemma}
\newcommand\bl{\begin{lemma}}
	\def\el{\end{lemma}}
\newtheorem{prop}[thn]{Proposition}
\newcommand\bp{\begin{prop}}
	\def\ep{\end{prop}}
\newtheorem{cor}[thn]{Corollary}
\newcommand\bc{\begin{cor}}
	\def\ec{\end{cor}}
\newtheorem{re}[thn]{Remark}
\newcommand\br{\begin{re}}
	\def\er{\end{re}}
\newtheorem*{defn}{Definition}
\newcommand\bd{\begin{defn}}
	\def\ed{\end{defn}}
\newtheorem{exam}[thn]{Example}
\newcommand\bxa{\begin{exam}}
	\def\exa{\end{exam}}

\def\be{\begin{equation}}
	\def\ee{\end{equation}}
\def\ba{\begin{array}}
	\def\ea{\end{array}}
\def\bdsc{\begin{description}}
	\def\edsc{\end{description}}
\def\bnum{\begin{enumerate}}
	\def\enum{\end{enumerate}}
\def\bea{\begin{eqnarray*}}
	\def\eea{\end{eqnarray*}}
\def\ben{\begin{eqnarray}}
	\def\een{\end{eqnarray}}
\def\btab{\begin{tabular}}
	\def\etab{\end{tabular}}
\def\bmat{\begin{bmatrix}}
	\def\emat{\end{bmatrix}}

\def\<{\langle}
\def\>{\rangle}

\def\ts{\textstyle}

\def\wh{\widehat }

\def\B{\mathbb{B}}
\def\Bq{{}^{\tt q}\kern-1.7pt\mathbb{B}}
\def\Be{\mathbb{E}}
\def\Dq{{}^{\tt q}\kern-2pt{D}}
\def\De{{}^{\epsilon}\kern-1pt{D}}
\def\Lq{{}^{\tt q}\kern-2.7pt\L}
\def\bA{\mathbb{A}}

\def\L{\mathfrak{L}}

\let\e\relax
\newcommand{\e}{\boldsymbol{e}}

\newcommand{\0}{\textbf{0}} 

\let\u\relax
\newcommand{\u}{\boldsymbol{u}}

\newcommand{\x}{\boldsymbol{x}}

\newcommand{\z}{\boldsymbol{z}}
\newcommand{\1}{\mathds{1}}

\newcommand{\calL}{\mathcal{L}}
\newcommand{\rvec}[1]{\boldsymbol{#1}}
\newcommand{\btau}{\boldsymbol{\tau}}
\newcommand{\btauq}{{}^{\tt q}\kern-2pt\boldsymbol{\tau}}

\newcommand{\qd}[1]{\boldsymbol{\{}#1\boldsymbol{\}}}

\begin{document}
	
\begin{frontmatter}		
\title{Inverse Formulae for $q$-analogues of Bipartite Distance Matrix}
\author[1]{R. Jana\corref{cor1}\fnref{rakesh}}
\ead{rjana@math.iitb.ac.in}
\address[1]{Department of Mathematics, IIT Bombay, India}
\cortext[cor1]{Corresponding author.}

\begin{abstract} 
We consider two distinct $q$-analogues of the bipartite distance matrix, namely the $q$-bipartite distance matrix and the exponential distance matrix. We provide formulae of the inverse for these matrices, which extend the existing results for the bipartite distance matrix. These investigations lead us to introduce a $q$-analogue version of the bipartite Laplacian matrix.
\end{abstract}
		
\end{frontmatter}
	

\section{Introduction}

The distance matrix finds numerous applications across various fields, including chemistry, molecular biology, and telecommunications. In 1971, Graham and Pollack \cite{GrPo} presented a remarkable result that the determinant of the distance matrix of a tree depends solely on the number of vertices in the tree. Specifically, for a tree $T$ with $n$ vertices, the determinant of its distance matrix is given by $\det D(T) = (-1)^{n-1}2^{n-2}(n-1)$. This groundbreaking discovery sparked significant interest among researchers.
In the same paper \cite{GrPo}, Graham and Pollack established a connection between the loop addressing problem in data communication systems and the count of negative eigenvalues of the distance matrix. Subsequently, in 1977, Graham, Hoffman, and Hosoya \cite{grahametall} demonstrated that the distance matrix $D(G)$ of a graph $G$ depends solely on its constituent blocks, regardless of how these blocks are connected. In 1978, Graham and Lov\'asz \cite{GrLo} provided a formula to compute the inverse of the distance matrix for trees.
Bapat \cite{bapat06} and Bapat et al. \cite{BKN} extended many results related to the distance matrix of trees to weighted trees. In 2006, Bapat, Lal, and Pati \cite{BPL} introduced two types of $q$-analogue versions of the distance matrix: the $q$-distance matrix and the exponential distance matrix. These $q$-analogue versions generated considerable interest and were explored by various researchers (see, for example, \cite{BPL, siva10, BS11, HLJ}). Let us revisit the definitions of the $q$-distance matrix and exponential distance matrix for a graph $G$.

Let $q$ be an indeterminate. For a positive integer $k$, we define $\qd{k} := 1 + q + q^2 + \ldots + q^{k-1}$. Let us set $\qd{0} = 0$. The square matrix $\Dq(G) := [\qd{\dist_G(i,j)}]$ is known as the {\em $q$-distance matrix} of $G$, and the square matrix $\De(G) := [q^{\dist_G(i,j)}]$ is the {\em exponential distance matrix} of $G$. Notably, ${}^1\kern-2ptD(G) = D(G)$. Hence, the $q$-distance matrix is a generalization of the usual distance matrix. 

The determinants of the $q$-distance matrix and exponential distance matrix are provided in \cite{BPL, YY07}.

\bt\label{th: det of Dq and De}
	Let $T$ be a tree with $n$ vertices. Let $\Dq$ and $\De$ be the $q$-distance matrix and the exponential distance matrix of $T$, respectively. Then $\det(\Dq) = (-1)^{n-1}(n-1)(1+q)^{n-2}$ and $\det(\De) = (1-q^2)^{n-1}$.
\et

For bipartite graphs, we can use a smaller matrix, known as the {\em bipartite adjacency matrix}, to store adjacency information efficiently (see \cite{BMbook}). Consider a labelled, connected, bipartite graph $G$ with the vertex bipartition $(L = \{l_1, \ldots, l_m\}, R = \{r_1, \ldots, r_n\})$. The \textit{bipartite adjacency matrix} $\bA(G)$ of $G$ is an $m \times n$ matrix, with rows indexed by $l_1, \ldots, l_m$ and columns indexed by $r_1, \ldots, r_n$. The $(i, j)$th entry of $\bA(G)$ is $1$ if $l_i \sim r_j$, and $0$ otherwise. The \textit{bipartite distance matrix} $\B(G)$ of $G$ is an $m \times n$ matrix, with rows indexed by $l_1, \ldots, l_m$ and columns indexed by $r_1, \ldots, r_n$. The $(i, j)$th entry of $\B(G)$ is $\dist(l_i, r_j)$.

Significant research has been conducted on the bipartite adjacency matrix of bipartite graphs with unique perfect matchings (see \cite{Godsil1985, PP15, PP15_2, PP17, Yang_Ye_2017} and references therein). The bipartite distance matrix for nonsingular trees, i.e., trees with perfect matchings, studied in \cite{bjp}. It was shown that unlike the determinant of the distance matrix of a tree, which is independent of the tree's structure, the determinant of the bipartite distance matrix of a nonsingular tree depends on its structure. A combinatorial expression for the determinant of the bipartite distance matrix of a nonsingular tree presented in  \cite{bjp}. Remarkably, the bipartite distance matrix of a nonsingular tree is always invertible, and its inverse is provided in  \cite{BJP2023}. During the investigation of the inverse of the bipartite distance matrix, the authors in \cite{BJP2023} uncovered a nontrivial generalization of the usual Laplacian matrix for a tree, which they have termed the {\em bipartite Laplacian matrix.} (We provide the definition of the bipartite Laplacian matrix in Section \ref{sec: inv-bd-matrix}.) Unlike the usual Laplacian matrix, the bipartite Laplacian matrix is not necessarily symmetric, but it possesses many elementary properties akin to the usual Laplacian matrix, see Theorem \ref{th: bl-properties}.

In \cite{Jana2023}, Jana studied of two distinct $q$-analogue versions of the bipartite distance matrix, namely, the $q$-bipartite distance matrix and the exponential bipartite distance matrix. Let's revisit the definitions of these matrices for a bipartite graph $G$.

Consider a connected, labelled bipartite graph $G$ with a vertex bipartition $(L = {l_1, \ldots, l_m}, R = {r_1, \ldots, r_n})$. The \textit{$q$-bipartite distance matrix} $\Bq(G)$ (respectively, the \textit{exponential bipartite distance matrix} $\Be(G)$) of $G$ is an $m \times n$ matrix, where the rows are indexed by $l_1, \ldots, l_m$, and the columns are indexed by $r_1, \ldots, r_n$. In $\Bq(G)$ (respectively, $\Be(G)$), the $(i, j)$th entry is $\qd{\dist_T(l_i,r_j)}$ (respectively, $q^{\dist_T(l_i,r_j)}$).

Let $G$ be a connected, labelled, bipartite graph with the vertex bipartition $(L = \{l_1, \ldots, l_m\}, R = \{r_1, \ldots, r_n\})$. The \textit{$q$-bipartite distance matrix} $\Bq(G)$ (respectively, \textit{exponential bipartite distance matrix} $\Be(G)$) of $G$ is an $m \times n$ matrix, with rows indexed by $l_1, \ldots, l_m$ and columns indexed by $r_1, \ldots, r_n$. The $(i, j)$th entry of $\Bq(G)$ (respectively, $\Be(G)$) is $\qd{\dist_T(l_i,r_j)}$ (respectively, $q^{\dist_T(l_i,r_j)}$).

If $q=1$ then $\Bq(G)=\B(G)$ and $\Be(G)=\1\1^t$. Therefore, the $q$-bipartite distance
matrix  is a generalization of the bipartite distance matrix. 

Let $T$ be a nonsingular tree on $2p$ vertices. It is shown in \cite{Jana2023} that $\det\Be(T)=q^p(1-q^2)^{p-1}$ and $\det\Bq(T)$ is divisible by $q^{p-1}(1+q)^{p-1}$.  Define $\bdq(T)$, the $q$-bipartite distance index of $T$, as $\bdq(T):=\frac{\det\Bq(T)}{q^{p-1}(1+q)^{p-1}}.$ Therefore, $\det\Bq(T)=q^{p-1}(1+q)^{p-1}\bdq(T)$. The $q$-bipartite distance index of $T$ follows an inclusion-exclusion type principle, enabling the recursive calculation of the $q$-bipartite distance index of any nonsingular tree. For more details, please refer to \cite[Theorem 20]{Jana2023}. Additionally, a standalone combinatorial formula for the $q$-bipartite distance index of a tree $T$ is provided in \cite[Theorem 23]{Jana2023}, by using the $f$-alternating sum technique (see \cite[Section 5]{bjp}).

It is evident that for a nonsingular tree $T$, the exponential bipartite distance matrix of $T$ is invertible when $q\neq 0,\pm 1$. Moreover, when $q\neq 0,- 1$ and $\bdq(T)\neq 0$, the $q$-bipartite distance matrix of $T$ is also invertible. In this article, our main focus is to provide the inverse formulas for both the exponential bipartite distance matrix and the $q$-bipartite distance matrix of a nonsingular tree whenever they exist.

We organize the article as follows.
In Section \ref{sec: inv-q-bd-matrix}, we discuss how to obtain the inverse of the bipartite distance matrix of a nonsingular tree using the bipartite Laplacian matrix.
In Section \ref{sec: q-bl-matrix}, we define the $q$-analogous version of the bipartite Laplacian matrix and explore its behaviour when extending a nonsingular tree by attaching a new $P_2$ at some vertex of the old tree.
In Section \ref{sec: inv-exp-bd}, we study the inverse of the exponential bipartite distance matrix of a nonsingular tree. In Section \ref{sec: inv-q-bd-matrix}, we introduce $q$-analogous versions of a few more vectors and finally provide the inverse formula for the $q$-bipartite distance matrix of a nonsingular tree.

\section{Inverse of the bipartite distance matrix}\label{sec: inv-bd-matrix}
Recently, in \cite{BJP2023}, it is shown that the inverse of the bipartite distance matrix of a nonsingular tree $T$ is a rank one update of a Laplacian matrix, which is referred to as the bipartite Laplacian matrix. Remarkably, the usual Laplacian matrix of any tree can be seen as a special case of the bipartite Laplacian matrix.
 Before providing the definition of the bipartite Laplacian matrix for a nonsingular tree, let's revisit some  definitions.

Let $T$ be a nonsingular tree on $2p$ vertices with a standard vertex bipartition $(L,R)$. An alternating path in $T$ is called an {\em even alternating path} if it contains even number of matching edges and an alternating path in $T$  is called an {\em odd alternating path} if it contains odd number of matching edges.

\bd\rm
Let $T$ be a nonsingular tree $T$ on $2p$ vertices and $(L,R)$ be a standard vertex 
bipartition of $T$. The {\em  bipartite Laplacian matrix} of  $T$, denoted by  
$\L(T)$ or simply by $\L$ is the $p\times p$ matrix whose rows are indexed by 
$r_1,\ldots, r_p$ and the columns are indexed by $l_1,\ldots, l_p$. The  
$(i,j)$th entry of $\L(T)$, denoted by $\L_{ij}$, is defined as 
$$
\L_{ij}=\left\{\ba{rl} 
d(r_i)d(l_i)-1 & \mbox{ if } i=j;\\
d(r_i)d(l_j) & \mbox{ if $i\neq j$ and the $r_i$-$l_j$ path is an odd alternating path};\\
-d(r_i)d(l_j) & \mbox{ if $i\neq j$ and the $r_i$-$l_j$ path is an even alternating path};\\
-1 & \mbox{ if $i\neq j$ and $r_i\sim l_j$};\\
0 & \mbox{ otherwise}.
\ea\right.
$$
\ed
Here we want to emphasize that in the definition of $\B(T),\Bq(T)$, and $\Be(T)$ we indexed rows by $l_1,\ldots,l_p$ and columns by $r_1,\ldots, r_p$. But in the definition of bipartite Laplacian matrix we indexed rows by $r_1,\ldots,r_p$ and columns by $l_1,\ldots, l_p$.

In the following result we highlight some similarities between the bipartite Laplacian matrix and the usual Laplacian matrix.

\bt\cite[Theorem 9]{BJP2023}\label{th: bl-properties} 
Let $T$ be a nonsingular tree on $2p$ vertices with a standard vertex bipartition $(L,R)$. 
Suppose $\L$ is the bipartite Laplacian matrix of $T$. Then the following assertions hold.
\begin{enumerate}[(a)]
	\item The row and the column sums of  $\L$ are zero. (A similar property also holds for
	the usual Laplacian matrix of any graph.)
	
	\item The cofactors of any two elements of $\L$ are equal to one.  (In the case of the usual Laplacian matrix of a graph, the cofactors of any two elements are equal to the number of spanning trees.)
	
	\item The rank of $\L$ is $p-1$. (A similar result holds for the usual Laplacian matrix of a connected graph.)
	
	\item  The algebraic multiplicity of $0$ as an eigenvalue of $\L$ is one.   (This property is also true for the usual Laplacian matrix of a connected graph.)
	
	\item  If $\u$ is an eigenvector of $\L$ corresponding to 
	an eigenvalue $\lambda\neq 0$  then $\1^t\u=0$.  (A similar property holds for the usual Laplacian matrix of any graph.)

	\item If $T=F\circ K_1$ for some tree $F$ on $p$ vertices. Then $\calL(F)=\L(T)$ where 
	$\calL$ is the usual Laplacian matrix of $F$. (This means that the usual Laplacian matrix of a tree $F$ can be seen as a bipartite Laplacian matrix of another tree $T$.)
	
	\item The matrix $\L$ is a symmetric matrix if and only if $T$ is a corona tree, that is,
	if $T=F\circ K_1$ for some tree $F$. 
	 (In contrast, the usual Laplacian matrix is always symmetric.)
\end{enumerate}
\et 

It's worth noting that in their work \cite{BJP2023}, Bapat, Jana, and Pati put forward a conjecture regarding the bipartite Laplacian matrix of a nonsingular tree. They proposed that this matrix is diagonalizable, and all of its eigenvalues are nonnegative real numbers. This conjecture, coupled with the properties outlined in Theorem \ref{th: bl-properties}, highlights the potential significance and relevance of the bipartite Laplacian matrix in the realm of nonsingular trees.

Similar to the inverse of the usual distance matrix of a tree, which can be viewed as a rank one update of its Laplacian matrix as shown in \cite{BKN}, the bipartite distance matrix of a nonsingular tree follows a similar pattern. It can also be seen as a rank one update of its bipartite Laplacian matrix. Before stating the result let us define some useful terminologies.

Let $v$ be a vertex in a nonsingular tree $T$. By $\c A_{T,v}^+$ we denote the set of all even alternating paths in $T$  that start at
$v$. Similarly, we denote the set of all odd alternating paths in $T$ that start at 
$v$ by $\c{A}_{T,v}^-$. By $\diff_T(v)$, we mean the quantity 
$\diff_T(v):=|\c{A}_{T,v}^+|-|\c{A}_{T,v}^-|$.   Let us
define the vector $\btau$ by $\btau(v):=1-d(v)(1+\diff(v))$ for each $v\in T$. 

\bt \label{th: inverse of B}
Let $T$ be a nonsingular tree on $2p$ vertices with a standard vertex bipartition $(L,R)$. 
Let $\B(T)$ and $\L(T)$ be the bipartite  distance matrix and the bipartite Laplacian matrix 
of $T$, respectively. 
Let $\btau_r(T)$ and $\btau_l(T)$ be the restriction of the vector $\btau(T)$ on $R$ and $L$, respectively.  
Then 
$$
\B(T)^{-1}=-\frac{1}{2}\L(T)+\frac{1}{\bdf(T)}\btau_r(T)\btau_l^t(T).
$$
\et

In this article, we mostly employ the proof strategy used in \cite{BJP2023}, which involves using information from the existing tree to establish the inductive case for a tree formed by attaching a $P_2$ to a vertex in the old tree. Below we provide the definition of attaching a $P_2$ to a vertex.

\bd[Attaching a  $P_2$ to a vertex]\rm (a)  Consider a tree $T$ and a vertex $v$. Let $\wh T$ denote the tree resulting from $T$ by adding two new vertices, $u$ and $w$, along with the edges $[v,u]$ and $[u,w]$. We refer to this operation as {\em attaching a new $P_2$ at the vertex $v$.}

(b) For a nonsingular tree $T$ with a vertex set of size $2p$ and a standard vertex bipartition $(L,R)$, let $\wh T$ be the tree formed by attaching a new $P_2$ to some vertex $v\in T$. To compute its $q$-bipartite distance matrix, exponential distance matrix and the bipartite Laplacian matrix, we label the new vertices according to the following procedure:

\bdsc
\item{i)} if $v\in L$, then we put $u=r_{p+1}$, $w=l_{p+1}$, and
\item {ii)} if $v\in R$, then we put $u=l_{p+1}$, $w=r_{p+1}$.
\edsc
\ed

\section{$q$-bipartite Laplacian Matrix}\label{sec: q-bl-matrix}
As our primary focus is to provide the inverse of the $q$-bipartite distance matrix, we begin by introducing the $q$-analogous version of the bipartite Laplacian matrix, which we call the {\em $q$-bipartite Laplacian matrix}. For a positive integer $k$, we define $k_q$ as $1+(k-1)q^2$.  
\bd\rm
Let $T$ be a nonsingular tree $T$ on $2p$ vertices and $(L,R)$ be a standard vertex 
bipartition of $T$. The {\em $q$-bipartite Laplacian matrix} of  $T$, denoted by  
$\Lq(T)$ or simply by $\Lq$ is the $p\times p$ matrix whose rows are indexed by 
$r_1,\ldots, r_p$ and the columns are indexed by $l_1,\ldots, l_p$. The  
$(i,j)$th entry of $\Lq(T)$, denoted by $\Lq_{ij}$, is defined as 
$$
\Lq_{ij}=\left\{\ba{rl} 
d(r_i)_qd(l_i)_q-q^2 & \mbox{ if } i=j;\\
d(r_i)_qd(l_j)_q & \mbox{ if $i\neq j$ and the $r_i$-$l_j$ path is an odd alternating path};\\
-d(r_i)_qd(l_j)_q & \mbox{ if $i\neq j$ and the $r_i$-$l_j$ path is an even alternating path};\\
-q^2 & \mbox{ if $i\neq j$ and $r_i\sim l_j$};\\
0 & \mbox{ otherwise}.
\ea\right.
$$
\ed

Clearly, for a nonsingular tree $T$, if we put $q=1$, then $\Lq=\L$, the bipartite Laplacian matrix of $T$. Therefore, the $q$-bipartite Laplacian matrix is a generalization of the usual bipartite Laplacian matrix.

Similar to signed degree vector as defined in \cite{BJP2023}, we now  introduce the concept of a $q$-analogues of the signed degree vector which relates the structure of the bipartite Laplacian 
of the new tree with that of the old one.  

\bd\rm  
Let $T$ be a nonsingular tree on $2p$ vertices with the standard vertex bipartition
$(L,R)$ and $v$ be a vertex. Then the {\em $q$-signed degree vector} $\qsdv{v}$ at $v$ is 
defined in the following way.
\bnum
\item If $v\in L$, then for $i=1,\ldots,p$, we define
\bdsc
\item{i)} $\qsdv{v}(i)=d_T(r_i))_q$ if the $v$-$r_i$ path is an odd alternating path, 
\item{ii)} $\qsdv{v}(i)=-(d_T(r_i))_q$ if the $v$-$r_i$ path is an even alternating path, and 
\item{iii)} $\qsdv{v}(i)=0$ if the $v$-$r_i$ path is not an alternating path.
\edsc

\item Similarly, if $v\in R$, then for  $i=1,\ldots,p$, we define
\bdsc
\item{i)} $\qsdv{v}(i)=(d_T(l_i))_q$ if the $v$-$l_i$ path is an odd alternating path, 
\item{ii)} $\qsdv{v}(i)=-(d_T(l_i))_q$ if the $v$-$l_i$ path is an even alternating path, and 
\item{iii)} $\qsdv{v}(i)=0$ if the $v$-$l_i$ path is not an alternating path.
\edsc 
\enum
\ed

Clearly, if $q=1$, then $\qsdv{v}=\sdv{v}$ for each $v\in T$.

In the following result, we discuss how the structure of $\Lq$ changes after attaching a $P_2$ to a vertex.

\bl\label{lem: add at l} 
Let $T$ be a nonsingular tree on $2p$ vertices with a standard vertex bipartition $(L,R)$.
Let $\widehat T$ be the tree obtained from $T$ by attaching a new $P_2$ at $v$.
Let $\qsdv{v}$ be the signed degree vector at $v$ of $T$. 
\bnum[(a)]
\item If $v=l_k$ for some $k$, then $
\Lq(\widehat T)= \begin{bmatrix}
	\Lq(T)+q^2\qsdv{v} \e_k^t & -\qsdv{v}\\
	-q^2\e_k^t & 1
\end{bmatrix}
$.

\item If $v=r_k$ for some $k$, then 
$
\Lq(\widehat T)= \begin{bmatrix}
	\Lq(T)+q^2\e_k\qsdv{v}^t & -q^2\e_k\\
	-\qsdv{v}^t & 1
\end{bmatrix}
$.
\enum
\el	

\begin{proof}
We only provide the proof of item (a) as the proof of item (b) can be dealt similarly.
Without loss of any generality, let us assume that $\widehat T$ obtained from $T$ by adding  a new path $[l_k, r_{p+1},l_{p+1}]$ for some $1\leq k\leq p$. Clearly $\wh L=L\cup\{l_{p+1}\}$ and $\wh R=R\cup\{r_{p+1}\}$ is a standard vertex bipartition of $\wh T$.  Let  $\Lq$ and $\wh \Lq$ be the bipartite Laplacian 
matrix of $T$ and $\widehat T$, respectively.  Since $[r_{p+1},l_{p+1}]$ is the only 
alternating path  that starts at $r_{p+1}$ and $d_{\widehat T}(r_{p+1})=2$ with $[r_{p+1}, 
l_k]$ is not a  matching edge, it follows that the all entries of the $(p+1)$th row of 
$\wh \Lq$ is  zero except $\wh \Lq(p+1,p+1)=d_{\widehat T}(r_{p+1})_qd_{\widehat T}(l_{p+1})_q-1=1$ 
and $\wh \Lq(p+1,k)=-1$. Hence $\wh \Lq(p+1,:)=\begin{bmatrix}-q^2e_k^t&1\end{bmatrix}$.

Let us take $i=1,\ldots,p$.  Then $r_i\nsim l_{p+1}$. Note that the $l_k$-$r_i$ path is an odd 
alternating path if and only if the $l_{p+1}$-$r_i$ path is an even alternating path. Similarly, the $l_k$-$r_i$ path is an even alternating path if and only if the  
$l_{p+1}$-$r_i$ path is an odd alternating path. Since $d_{\widehat T}(l_{p+1})=1$, it follows that $\wh \Lq(\{1,\ldots,p\},p+1)=-\qsdv{l_k}$, where $\qsdv{l_k}$ is the $q$-signed degree vector of $T$ at $l_k$.

Since $d_T(u)=d_{\widehat T}(u)$ for each $u\in T$ other than $l_k$, it follows that $\wh \Lq(i,j)=\Lq(i,j)$ for each $i=1,\ldots,p$ and $j=1,\ldots,k-1,k+1,\ldots,p$. 

Finally, notice that  $d_{\widehat T}(l_k)=d_T(l_k)+1$.  Therefore, for $i=1,\ldots,p$, we have
$$
\wh\Lq(i,k)= \left\{\ba{rl} 
d_{T}(r_i)_q(d_T(l_k)_q+q^2)-q^2 & \mbox{ if } i=k;\\
d_T(r_i)_q(d_T(l_k)_q+q^2) & \mbox{ if $i\neq k$ and the $r_i$-$l_k$ path is an odd alternating path};\\
-d_T(r_i)_q(d_T(l_k)_q+q^2) & \mbox{ if $i\neq k$ and the $r_i$-$l_k$ path is an even alternating path};\\
-q^2 & \mbox{ if $i\neq k$ and $r_i\sim l_k$};\\
0 & \mbox{ otherwise}.
\ea\right.
$$ 
Therefore, $\wh\Lq(\{1,\ldots,p\},k)=\Lq(\{1,\ldots,p\},k)+q^2\qsdv{l_k}$. This completes the proof.
\end{proof}

In the next result we extend the result \cite[Lemma 7]{BJP2023} which tells that the sum of all entries in a signed degree vector is always one.

\bl\label{lem: sum of mu}
Let $T$ be a nonsingular tree on $2p$ vertices with a standard bipartition $(L,R)$.
Let $u$ be any vertex in $T$ and $\sdv{u}$ be the signed degree vector at $u$. Then
$\1^t\sdv{u}=(\diff(u)+1)q^2-\diff(u)$.   
\el  

\begin{proof}
	We proceed by induction on $p\geq 1$.
	For $p=1$ the result is trivial. Assume the result to be true for nonsingular trees with
	less than $2p$ vertices. Let $T$ be a nonsingular tree  on $2p$ vertices with a 
	standard bipartition $(L,R)$. Let $u\in R$. (The case of $u\in L$ can be dealt 
	similarly.) Let $\qsdv{}$ be the $q-$signed degree vector of $u$ in $T$.
	
	Suppose $[v_0,v_1,\ldots,v_k]$ is a longest path in $T$. As $p>1$, we have $k\ge 3$ and
	so we may assume that $v_0, v_1\neq u$. As $T$ is nonsingular and this is a longest 
	path, we have $d(v_0)=1$ and $d(v_1)=2$. Without loss of any  generality, let us assume  
	$v_0,v_1\in\{l_p,r_p\}$. Let $\wh T=T-\{v_0,v_1\}$ be the tree obtained from $T$ 
	by removing the vertices $v_0$ and $v_1$. Clearly, $u\in \wh T$. Let $\wh{\qsdv{}}$ be the signed 
	degree vector of $u$ in $\wh T$. Note that $\wh{\qsdv{}}$ is vector of size $p-1$. Clearly, 
	$d_T(v)=d_{\wh T}(v)$ for each $v\in\wh T-v_2$ and $d_T(v_2)=d_{\wh T}(v_2)+1$. 
	It follows that $\qsdv{}(i)=\wh{\qsdv{}}(i)$ for each $l_i\in L\setminus\{v_2\}$.
	
	If either $v_2\in R$ or the $u$-$v_2$ path is not an alternating path then 
	$\qsdv{}=\begin{bmatrix}\wh{\qsdv{}}& 0\end{bmatrix}$ and the  result follows by induction. 
	Now we assume that $v_2\in L$ and the  $u$-$v_2$ path is an alternating path. Then 
	$v_2\sim r_p$ and $d_T(l_p)=1$. Let $v_2=l_k$ for some $1\leq k<p$. Note that the $u$-$l_p$ 
	path is also an alternating path and so we have $\wh{\rvec{\mu}}(k)=(-1)^{t} d_{\wh T}(v_2)_q$ 
	for some $t$  and $\rvec \mu(p)=(-1)^{t+1}$. Since  $\rvec\mu(k)=(-1)^td_T(v_2)_q=\wh{\rvec{\mu}}(k)+(-1)^tq^2$, 
	it follows that 
	$$\rvec\mu=\begin{bmatrix}\wh{\rvec{\mu}}^t& (-1)^{t+1}\end{bmatrix}^t+(-1)^{t}\begin{bmatrix}q^2e_k& 0\end{bmatrix}^t.$$ 
	Clearly, $\ts \diff_{T}=\diff_{\wh T}+(-1)^t$. 
	Hence, the result follows by applying  induction. 
\end{proof}

\bt 
Let $T$ be a nonsingular tree on $2p$ vertices and $\Lq$ be the $q$-Laplacian matrix. Then $\det\Lq=1-q^2$
\et 

\begin{proof}
	We use induction on $p$. The result clearly true for $n=1$ and $2$. Let us assume the result to be true for nonsingular tree on $2p$ vertices. Let $\wh T$ be a nonsingular tree on $2p+2$ vertices. Then $\widehat T$ can be viewed as obtained from some nonsingular tree $T$ with $2p$ vertices by attaching a new $P_2$ at a vertex $v$. Without loss of generality, let's assume that $v=l_k$ for some $1\leq k\leq p$.
	Let $\qsdv{}$ be the $q$-signed degree vector at $l_k$ of $T$. By Lemma \ref{lem: add at l}, we have
	$$\Lq(\widehat T)= 
	\bbm 
	\Lq(T)+q^2\qsdv{v} \e_k^t & -\qsdv{v}\\
	-q^2\e_k^t & 1
	\ebm
	$$
	By using Schur's formula for the determinant, we get 
	$$
	\det\Lq(T)=\det\big( \Lq(T)+q^2\qsdv{v} \e_k^t -   q^2 \qsdv{v} \e_k^t \big) = \det\Lq(T).
	$$
	Hence, the result follows by applying induction hypothesis. 
\end{proof}

In the following remark we discuss how the $q$-bipartite Laplacian matrix of a nonsingular
tree can be obtained from some of its nonsingular subtrees. This plays a crucial rule proving our subsequent results.

\br \label{remark: understanding}\rm 
Consider the tree $T$ with a matching edge $[l_{k_1},r_{k_1}]$, see Figure \ref{fig: understanding}. 
Let the degree of $l_{k_1}$ be $s$, $s\geq 1$. Let  $r_{{k_1}+1},r_{k_2+1},\ldots,r_{k_{s-1}+1}$ be
some distinct vertices other than $r_{k_1}$ that are adjacent to $l_{k_1}$. 
Note that when we delete the edges $[l_{k_1},r_{k_1+1}]$, $\ldots,  [l_{k_1},r_{k_{s-1}+1}]$, 
we obtain $s$ many smaller nonsingular trees, say $T_1,\ldots, T_s$. 
Assume that the vertex set of $T_1$  is $\{l_1,\ldots,l_{k_1},r_1,\ldots,r_{k_1}\}$, the vertex set of 
$T_2$ is $\{l_{k_1+1},\ldots,l_{k_2},r_{k_1+1},\ldots,r_{k_2}\}$, and so on up to  
the vertex set of $T_s$  is $\{l_{k_{s-1}+1},\ldots,l_{k_s},r_{k_{s-1}+1},\ldots,r_{k_s}\}$. 
Let us put an arrow on the edge $[l_{k_1},r_{k_1}]$ from $r_{k_1}$ to $l_{k_1}$. 
This arrow indicates that, from a vertex $r_i$ in $T_2$, we do not have an alternating
path to a vertex in $T_1$. Similarly,  from a vertex $r_i$ in $T_3$, we do not have an
alternating path to a vertex in $T_1,T_2,T_4,T_5, \ldots, T_s$. 
Similar statements are true for vertices $r_i$ in $T_4, \ldots, T_s$. 
Also, from a vertex $l_i$ in $T_1$, we only have alternating paths to vertices in $T_1$ 
but not to a vertex in $T_2,\ldots,T_s$. Let us take $F_1$ be the tree $T_1$. 
For $i=2,\ldots, s$, let $F_i$ be the subtree of $T$ obtained by taking $F_{i-1}$ 
and $T_i$ and by inserting the edge $[l_{k_1}, r_{k_{i-1}+1}]$. 
Clearly $F_s$ is the original tree $T$.

\begin{figure}[h]
	\bpsp(0,6.75)\rput(3,3.5){\psline[showpoints=true,linestyle=dashed](3,0)(4,0)
		\psline[showpoints=true](4,0)(6,2)
		\psline[showpoints=true](6,-2.5)(4,0)(6,0)
		\psellipse(9,0)(3.2,.8)
		\psellipse(9,-2.5)(3.2,.8)
		\psellipse(9,2)(3.2,.8)
		\psellipse(1.5,0)(3,1)
		\rput(9,2.2){$l_{k_1+2},\ldots,l_{k_2}$}
		\rput(9,1.7){$r_{k_1+1},\ldots,r_{k_2}$}
		\rput(6.6,2){$r_{k_1+1}$}
		\rput(9,0.2){$l_{k_2+2},\ldots,l_{k_3}$}
		\rput(9,-.3){$r_{k_2+1},\ldots,r_{k_3}$}
		\rput(6.6,0){$r_{k_2+1}$}
		\rput(9.5,-2.3){$l_{k_{s-1}+2},\ldots,l_{k_s}$}
		\rput(9.5,-2.8){$r_{k_{s-1}+1},\ldots,r_{k_s}$}
		\rput(6.9,-2.5){$r_{k_{s-1}+1}$}
		\rput(0,.2){$l_1,\ldots,l_{k_1-1}$}
		\rput(0,-.3){$r_1,\ldots,r_{k_1-1}$}
		\rput(3,-.3){$r_{k_1}$}
		\rput(3.8,.3){$l_{k_1}$}
		\rput(1.5,.7){$T_1$}
		\rput(11.8,2){$T_2$}
		\rput(11.8,0){$T_3$}
		\rput(11.8,-2.5){$T_s$}
		\rput(5.8,-1){$\vdots$}
		\rput(1.8,-2.5){{\Large$T$}}
		\rput(3.5,0){\orient{0}}
		\psline[linestyle=dotted](-1,-2)(5,-.5)(5.3,1)(12.5,1)(12.5,3)(1,3)(-1.8,0)(-1,-2)
		\rput(3.3,2){$F_2$}
	}\epsp
	\caption{Understanding $\L(T)$.\label{fig: understanding}}
\end{figure}

a) Let $\qsdv{l_{k_1}}$ be the $q$-signed degree vector at $l_{k_1}$ of $T_1$ and let $\qsdv{}$ be the signed degree vector at $l_{k_1}$ of $T$. Then $\rvec\mu= \bbm \sdv{l_{k_1}} & \0\ebm^t$.

b) Let $\Lq(F_i)$ be the $q$-bipartite Laplacian matrix of $F_i$ for $i=1,\ldots,s$ and  $\Lq(T_i)$ be the $q$-bipartite Laplacian matrix of $T_i$ for $i=1,\ldots,s$. Clearly, $\Lq(T_1)=\Lq(F_1)$. Let $\qsdv{r_{k_i+1}}$ be the $q$-signed degree vector at $r_{k_i+1}$ of $T_{i+1}$, for $i=1,\ldots,s-1$. By $\rvec{E}^{ij}$ we denote the matrix of an appropriate size with $1$ at position $(i,j)$ and zero elsewhere. Then, for $i=2,\ldots,s$, we have 
$$\renewcommand{\arraystretch}{1.5}
\L(F_i)= \left[\ba{c|c|c|c} 
\L(T_1)+(i-1)q^2\sdv{l_{k_1}}\e_{k_1}^t & -\qsdv{l_{k_1}}\qsdv{r_{k_{1}+1}}^t &\cdots & 
-\qsdv{l_{k_1}} \qsdv{r_{k_{i-1}+1}}^t \\[.5em]\hline
-q^2\rvec{E}^{1k_1} & \L(T_2)+q^2\e_1\qsdv{r_{k_1+1}}^t & \cdots & \0 \\[.5em]\hline
\vdots &\0 &\ddots & \0\\\hline 
-q^2\rvec{E}^{1k_1} & \0 & \cdots & \L(T_i)+q^2\e_1\qsdv{r_{k_{i-1}+1}}^t 
\ea \right]
$$ 
In particular,
$$\renewcommand{\arraystretch}{1.5}
\L(T)=  \left[\ba{c|c|c|c} 
\L(T_1)+(s-1)q^2\sdv{l_{k_1}}\e_{k_1}^t & -\qsdv{l_{k_1}}\qsdv{r_{k_{1}+1}}^t &\cdots & 
-\qsdv{l_{k_1}} \qsdv{r_{k_{i-1}+1}}^t \\[.5em]\hline
-q^2\rvec{E}^{1k_1} & \L(T_2)+q^2\e_1\qsdv{r_{k_1+1}}^t & \cdots & \0 \\[.5em]\hline
\vdots &\0 &\ddots & \0\\\hline 
-q^2\rvec{E}^{1k_1} & \0 & \cdots & \L(T_s)+q^2\e_1\qsdv{r_{k_{s-1}+1}}^t 
\ea \right]
$$ 
\er

\section{Exponential bipartite distance matrix}	\label{sec: inv-exp-bd}

Let us recall the following result which tells that the exponential bipartite distance matrix of a nonsingular tree is independent of tree structure.

\bt\label{th: det of expB} \cite{Jana2023}
Let $T$ be a nonsingular tree on $2p$ vertices with a standard vertex bipartition $(L, R)$. Then $\det \Be(T)=q^{p}(1-q^2)^{p-1}$.
\et

By Theorem \ref{th: det of expB}, we observe that $\Be(T)$ is invertible whenever $q\neq 0,\pm 1$. In the following result, we present the inverse of $\Be(T)$ under the condition that $q\neq 0,\pm 1$.

\bt 
Let $T$ be a nonsingular tree on $2p$ vertices with a standard vertex bipartition $(L, R)$. Suppose $q\neq 0,\pm 1$. Then 
$$
\Be(T)^{-1}= \frac{\Lq}{q(1-q^2)}.
$$
\et

\begin{proof}
We proceed by induction on $p$. The base can be verified easily. 
 Assuming the result holds for $p$. Let $\widehat T$ be a nonsingular tree with $2p+2$ vertices. We can obtain $\widehat T$ from some nonsingular tree $T$ with $2p$ vertices by attaching a new $P_2$ at a vertex $v$. Without loss of generality, let's assume that $v=l_k$ for some $1\leq k\leq p$.
Let $\qsdv{}$ be the $q$-signed degree vector at $l_k$ of $T$. 
 
By item (a) of Lemma \ref{lem: add at l}, we have 
$$\Lq(\widehat T)= \begin{bmatrix}
	\Lq(T)+q^2\qsdv{} \e_k^t & -\qsdv{} \\
	-q^2\e_k^t & 1
\end{bmatrix}$$
Let $\x$ be a vector of size $p$ such that $\x(i)= q^{\dist(l_i,r_{p+1})}$ 
for each $i=1,\cdots, p$. 
Then the exponential bipartite distance matrix of $\widehat T$ can be written as
$$
\Be(\widehat T)= \bbm \Be(T) & \x\\
q^2\e_k^t\Be(T)&q\ebm.
$$
Now note that 
\begin{align}
	\Lq(\widehat T)\Be(\widehat T)&=  
	\bbm
		\Lq(T)+q^2\qsdv{} \e_k^t & -\qsdv{} \\
		-q^2\e_k^t & 1
	\ebm
	\bbm \Be(T) & \x\\
	q^2\e_k^t\Be(T)&q
	\ebm
	 \nonumber \\
	&=\bbm 
	\Lq(T)\Be(T) & \Lq(T)\x +q^2\qsdv{}\e_k^t\x-q\qsdv{}\\
	\0 & -q^2\e_k^t\x+q
	\ebm\nonumber  \\
	&=\bbm 
	q(1-q^2)I & \Lq(T)\x +q(q^2-1)\qsdv{}\\
	\0 & -q^2\e_k^t\x+q
	\ebm.\label{eq1}
\end{align}
The last equality follows from induction hypothesis and by using  $\e_k^t\x=\x(k)=\dist(l_k,r_{p+1})=1$.  Therefore, to complete the proof, we only need to show that   $\Lq(T)\x=q(1-q^2)\qsdv{}$.

Let the degree of $l_k$ in $T$ be $s$, ($s\geq 1$).
Let $T_1$ be the tree obtained from $T$ by removing all vertices adjacent to $l_k$ except the vertex $r_k$. 
If $s=1$ then $T_1$ is the same as $T$. 
Without loss of any generality, let us assume that $T_1$ has the vertex set $\{l_1,r_1,\ldots,l_{k-1},r_{k-1},l_k,r_k\}$.  
Let $\widehat{\qsdv{}}$ be the $q$-signed degree vector at $v$ of $T_1$.
By Remark \ref{remark: understanding}, $\qsdv{}=\bbm \widehat{\qsdv{}} & \0\ebm^t$.  
Further note that
\be\label{eq2}
\x=\Be(T)\e_k+(q^2-1)\matr{rrrrrrr}{q^{\dist(l_1,r_k)}&\cdots&q^{\dist(l_{k-1},r_{k})}&0&0&\cdots& 0}^t.
\ee 
Let $\z$ be a vector of size $(p-k)$ defined as follows. 
For $i=1,\ldots , (p-k)$,  $\z(i)=-1$ if $r_{k+i}$ adjacent to $l_k$ and $\z(i)=0$ otherwise.   

Hence, by Remark \ref{remark: understanding}, we have
\begin{align}
	\Lq(T)\bbm \Be(T_1)\e_k-q\e_k\\\0 \ebm^t  
	&= \bbm 
	\Lq(T_1)+(s-1)\widehat{\qsdv{}} \e_k^t & *\\
	q^2\z \e_k^t & * 
	\ebm \bbm \Be(T_1)\e_k-q\e_k\\\0 \ebm^t \nonumber\\
	& =\bbm 
	\Lq(T_1)\Be(T_1)\e_k + (s-1)\widehat{\qsdv{}}\e_k^t\Be(T_1)\e_k\\ q^3\z
	\ebm -
	\bbm 
	q\Lq(T_1)\e_k + q(s-1)\widehat{\qsdv{}}\\ q^3\z
	\ebm
	 \nonumber\\
	&= \bbm 
	q(1-q^2)\e_k -q(\widehat{\qsdv{}}-q^2\e_k) \\ \0
	\ebm \nonumber \\
	&=q(e_k-\qsdv{}). \label{eq3}
\end{align}
The second last equality holds by using the induction hypothesis and the fact that $\Lq(T_1)\e_k = \widehat{\qsdv{}}-q^2\e_k$. The last equality holds by using Remark \ref{remark: understanding}.
Hence, by \eqref{eq2} and \eqref{eq3}, it follows that 
$$
\Lq(T)\x = q(1-q^2) \e_k - q(1-q^2)(e_k-\qsdv{}) = q(1-q^2)\qsdv{}.
$$ 
This completes the proof.
\end{proof}

\section{$q$-bipartite distance matrix}	\label{sec: inv-q-bd-matrix}
In this section, first, we recall some terminologies from \cite{Jana2023} and using them we 
present a formula for the inverse of the $q$-bipartite distance matrix of a nonsingular tree.  

Let $T$ be a nonsingular tree on $2p$ vertices with a standard vertex bipartition $(L,R)$. 
The vector $\btauq (T)$, or simply $\btauq$, of size $2p$ is defined by 
$\btauq(v):=\big(1-d_T(v)\big)\big(1+\diff_T(v)\big)q^2-\diff_T(v)$ for each $v$ in $T$. The entries of $\btauq$ are 
ordered according to $l_1,\ldots,l_p,r_1,\ldots,r_p$. Clearly, for $q=1$, $\btauq$ is the vector $\btau$, as defined in \cite{bjp}.

By $\btauq_r(T)$, or simply by $\btauq_r$, we mean the restriction of $\btauq(T)$ on $R$. 
Similarly, by $\btauq_l(T)$, or simply by $\btauq_{l}$, we mean the restriction of 
$\btauq(T)$ on $L$.

The next result relates the structures the $\btauq_r$ vectors 
of the new tree with that of the old one under attaching a new $P_2$ at a vertex.
The first item was proved in \cite[Lemma 13]{Jana2023} and the proof of the second item is routine.
\bl\label{lem: prop of taur} 
Let $T$ be a nonsingular tree on $2p$ vertices with a standard vertex bipartition $(L,R)$.
Let $\wh T$ be the tree obtained from $T$ by attaching a new $P_2$ at $v$. 

(a) If $v=r_k$ for some $k$, then $\btauq_r(\wh T)=\bbm 
\btauq_r(T)\\0\ebm-(1+\diff_T(v))\bbm \e_k\\-1 \ebm$.

(b) If $v=l_k$ for some $k$, then $\btauq_r(\wh T)=\bbm \btauq_r(T)\\ 1\ebm - \bbm 
\qsdv{v}(T)\\0 \ebm$, where $\qsdv{v}(T)$ is the $q$-signed degree vector at $v$ of $T$. 
\el

In the following result, we present details regarding the row sums and column sums of the $q$-bipartite Laplacian matrix.

\bt\label{th: properties} 
Let $T$ be a nonsingular tree on $2p$ vertices with a standard vertex bipartition $(L,R)$. 
Suppose $\Lq(T)$ is the bipartite Laplacian matrix of $T$. Then we have
$$\1^t\Lq(T)=(1-q^2) (\btauq_l(T))^t\qquad  \text{and} \qquad \Lq(T) \1=(1-q^2) \btauq_r(T).$$
\et 

\begin{proof}
	We use induction on $p$. Let $\widehat T$ be a nonsingular tree on $2p+2$ vertices. Suppose $\widehat T$ is obtained from some nonsingular tree $T$ on $2p$ vertices by attaching a new 
	$P_2$ at some vertex $v$. Notice that either $v\in L$ or $v\in R$. (We shall proof the case  $v\in L$ as the other case can be dealt similarly.)
	
	Suppose $v=l_k$ and $l_k\sim r_{p+1}$. Clearly, $l_{p+1},r_{p+1}\notin T$. Let $\qsdv{}$ be the $q$-signed degree vector at $l_k$ of $T$. By Lemma \ref{lem: add at l}, we get 
	$$
	\Lq(\wh T)=
	\bbm 
	\Lq(T)+q^2\qsdv{} \e_k^t & -\qsdv{}\\
	-q^2\e_k^t & 1
	\ebm .
	$$
	Suppose $(z_l(\widehat T))_i = (1-q^2) \btauq(\wh T) (l_i)$ and  $(z_r(\widehat T))_i = (1-q^2) \btauq(\wh T) (r_i)$ for $1\leq i\leq p$.
	Clearly, $\ts \z_{\widehat T}(r_{p+1})= 1-q^2$ and $\z_{\widehat T}(l_{p+1})=\diff_{\widehat T}(l_{p+1})(q^2-1)$. By Lemma \ref{lem: sum of mu}, $1-\1^t\qsdv{} =(1+\diff_T(l_k))- (1+\diff_T(l_k))q^2 = \diff_{\widehat T}(l_{p+1})(q^2-1)$, as  $\ts \diff_{\widehat T}(l_{p+1})=-(1+\diff_T(l_k))$. Therefore, it follows that
	
	$$(\1^t\Lq(\wh T))_{p+1}= 1-\1^t\qsdv{}= (\z_l(\widehat T))_{p+1}\quad \text{and }\quad (\Lq(\wh T)\1)_{p+1}=1-q^2=(\z_r(\widehat T))_{p+1}.$$
	
	For each $1\leq i\leq p$, we have $d_{\wh T}(r_i)=d_{T}(r_i)$ and  $\diff_{\wh T}(r_i)=\diff_{T}(r_i)+\x(r_i)$, where $\x$ is a vector such that $\x(r_i)=0$ if the $r_i$-$l_{p+1}$ path is not an alternating path, $1$ if the $r_i$-$l_{p+1}$ path is an even alternating path, and $-1$ if the $r_i$-$l_{p+1}$ path is an odd alternating path. Therefore, by applying induction hypothesis, we get
	\begin{align*}
		\left(\Lq(T)\1 +(q^2-1)\qsdv{}\right)_i &= \ts 
		(\diff_T(r_i)+1)(d_{T}(r_i)-1)q^4+\left[\diff_T(r_i)(2-d_T(r_i))+(1-d_T(r_i))\right]q^2\\
		& \ts \qquad -\diff_T(r_i) + (q^2-1)(1+(d_T(r_i)-1)q^2)\x(i)\\
		&= (1-q^2) \btauq(\wh T) (r_i). 
	\end{align*}
	
	Further, for  $1\leq i\leq p$,  note that $\diff_{\wh T}(l_i)=\diff_{T}(l_i)$,  $d_{\wh T}(l_i)=d_{T}(l_i)$ for $i\neq k$ and  $d_{\wh T}(l_k)=d_{T}(l_k)+1$.  Then, $(z_l(\widehat T))_i=(z_l(T))_i$ for each $1\leq i\leq p$ and $i\neq k$. Therefore, by induction hypotheses, it follows that $(\1^t\Lq(\wh T))_{i}=(z_l(\widehat T))_i$ for each $i\neq k$. For $i=k$, by using induction hypothesis and Lemma \ref{lem: sum of mu}, we have 
	\begin{align*}
		\left(\1^t\Lq(\widehat T)\right)_k &=  (\z_l(T))_k +q^2(\1^t\qsdv{}-1)\\
		&=\ts (\z_l(T))_k +q^2(q^2-1)(1+\diff_T(l_k)) = (1-q^2) \btauq(\wh T) (l_k)
	\end{align*}
	This completes the proof.
\end{proof}

\bt\cite{Jana2023}\label{th: main-result-Btau-bd} 
Let $T$ be a nonsingular tree on $2p$ vertices with a standard vertex bipartition $(L, R)$. Then following assertions hold. 

(a) $\Bq(T)\btauq_r(T)=\bdq(T)\1$ and $\btauq_l(T)^t\Bq(T)=\bdq(T)\1^t$. 

(b) Let $v$ be a vertex and let $\wh T$ be the tree obtained from $T$ by attaching a new $P_2$ at $v$. Then $\bdq(\wh T)= \bdq(T)+(1+q)(1+\ts\diff_T(v))$.
\et 

In the next result we discuss a relationship between the $q$-bipartite distance 
matrix and the $q$-bipartite Laplacian matrix of a nonsingular tree. 

\bl\label{lem: 111} 
Let $T$ be a nonsingular tree on $2p$ vertices with a standard vertex bipartition $(L,R)$.
Let $\Bq(T)$ and $\Lq(T)$ be the bipartite distance matrix and the bipartite Laplacian matrix 
of $T$, respectively. 
Let $\btauq_r(T)$ be the restriction of $\btauq(T)$ on  $R$.  Then 
$$
-\Lq(T)\Bq(T)+(1+q)\btauq_r(T)\1^t=q(1+q)I
$$
\el 

\begin{proof}
We proceed by induction on $p$. The base can be verified easily. 
Assuming the result holds for $p$. Let $\widehat T$ be a nonsingular tree with $2p+2$ vertices. We can obtain $\widehat T$ from some nonsingular tree $T$ with $2p$ vertices by attaching a new $P_2$ at a vertex $v$. Without loss of generality, let's assume that $v=l_k$ for some $1\leq k\leq p$.
Let $\qsdv{}$ be the $q$-signed degree vector at $l_k$ of $T$. 

By item (a) of Lemma \ref{lem: add at l}, we have 
$$\Lq(\widehat T)= \begin{bmatrix}
	\Lq(T)+q^2\qsdv{} \e_k^t & -\qsdv{} \\
	-q^2\e_k^t & 1
\end{bmatrix}$$
Let $\x$ be a vector of size $p$ such that $\x(i)= \qd{\dist(l_i,r_{p+1})}$ 
for each $i=1,\cdots, p$. 
Then the $q$-bipartite distance matrix of $\widehat T$ can be written as
$$
\Bq(\widehat T)= \bbm \Bq(T) & \x\\
(1+q)\1^t+q^2\e_k^t\Bq(T)&1\ebm.
$$
Now note that 
\begin{align}
	\Lq(\widehat T)\Bq(\widehat T)&=  
	\bbm
	\Lq(T)+q^2\qsdv{} \e_k^t & -\qsdv{} \\
	-q^2\e_k^t & 1
	\ebm
	\bbm \Bq(T) & \x\\
	(1+q)\1^t+q^2\e_k^t\Bq(T)&1
	\ebm
	\nonumber \\
	&=\bbm 
	\Lq(T)\Bq(T) -(1+q)\qsdv{}\1^t & \Lq(T)\x + (q^2-1)\qsdv{}\\
	(1+q)\1^t & 1-q^2
	\ebm \nonumber \\
	&= \bbm 
	(1+q)\btauq_r(T)\1^t-q(1+q)I -(1+q)\qsdv{}\1^t & \Lq(T)\x + (q^2-1)\qsdv{}\\
	(1+q)\1^t & 1-q^2
	\ebm 
	 \label{eq: eq1}
\end{align}
The last equality follows from induction hypothesis.  By part (b) of Lemma \ref{lem: prop of taur}, we get 
\begin{align*}
-\Lq(\wh T)\Bq(\wh T)+(1+q)\btauq_r(\wh T)\1^t & = -\Lq(\wh T)\Bq(\wh T)+(1+q)\bbm(\btauq_r(T)-\qsdv{})\1^t \\ \1^t \ebm \\
&= \bbm 
q(1+q)I & -\Lq(T)\x + (1+q)\btauq_r(T) -q(1+q)\qsdv{}\\
0 & q(1+q)
\ebm 
\end{align*}
Therefore, we only remain to show that  $\Lq(T)\x=(1+q)\btauq_r(T) -q(1+q)\qsdv{}$.

Let the degree of $l_k$ in $T$ be $s$, ($s\geq 1$).
Let $T_1$ be the tree obtained from $T$ by removing all vertices adjacent to $l_k$ except the vertex $r_k$. 
If $s=1$ then $T_1$ is the same as $T$. 
Without loss of any generality, let us assume that $T_1$ has the vertex set $\{l_1,r_1,\ldots,l_{k-1},r_{k-1},l_k,r_k\}$.  
Further, note that
\be\label{eq: eq2}
\x=\Bq(T)\e_k+ (1+q)\matr{rrrrrrr}{q^{\dist(l_1,r_k)}&\cdots&q^{\dist(l_{k-1},r_{k})}&0&0&\cdots& 0}^t.
\ee 

By \eqref{eq3}, we already have
\begin{align}
	\Lq(T)\bbm \Be(T_1)\e_k-q\e_k\\\0 \ebm^t  =q(\e_k-\qsdv{}). \label{eq: eq3}
\end{align}
By using induction hypothesis, $\Lq(T)\Bq(T)=(1+q)\btauq_r(T)\1^t-q(1+q)I$. 
Hence, by \eqref{eq: eq2} and \eqref{eq: eq3}, it follows that 
$$
\Lq(T)\x = (1+q)\btauq_r(T)-q(1+q)\e_k +q(1+q)(\e_k-\qsdv{})= (1+q)\btauq_r(T)-q(1+q)\qsdv{}.
$$ 
This completes the proof.
\end{proof}

We are now in a position to supply a formula for the inverse of the $q$-bipartite 
distance matrix of a nonsingular tree $T$.  
\bt 
Let $T$ be a nonsingular tree on $2p$ vertices with a standard vertex bipartition $(L, R)$. Suppose $q\neq 0,-1$ and $\bdq(T)\neq 0$. Then 
$$
\Bq(T)^{-1}= -\frac{1}{q(1+q)} \Lq(T) + \frac{1}{q\bdq(T)} \btauq_r(T)\btauq_l^t(T).
$$
\et

\begin{proof}
	By using Theorem \ref{th: main-result-Btau-bd}  and Lemma \ref{lem: 111} we have
	\begin{align*}
		\left(-\frac{1}{q(1+q)}\Lq(T)+\frac{1}{q\bdq(T)} \btauq_r(T)\btauq_l^t(T)\right)\Bq(T) = \frac{1}{q(1+q)}\left(-\Lq(T)\Bq(T) +(1+q)\btauq_r(T)\1^t\right) =I
	\end{align*}
	This completes the proof.
\end{proof}

\bibliographystyle{apalike}%
\bibliography{bibdatabase}

\end{document}